\documentclass[preprint,12pt,authoryear]{elsarticle}

\usepackage{amssymb}
\usepackage{amsmath}

\RequirePackage{fix-cm}
\RequirePackage{rotating}
\usepackage{booktabs} 
\usepackage{graphicx} 
\usepackage{xcolor} 
\usepackage{geometry} 
\usepackage{graphicx}

\usepackage{graphicx,epsfig}
\usepackage{amsthm}

\newtheorem{theorem}{Theorem} 
\newtheorem{corollary}{Corollary}

\usepackage{exscale}
\usepackage{booktabs}
\usepackage{algorithm}
\usepackage{algorithmicx}
\usepackage{algpseudocode}

\usepackage{multirow}
\usepackage{longtable}

\usepackage{rotating, tabularx}
\usepackage{enumitem}
\DeclareMathOperator*{\argmin}{arg\,min}

\newtheorem{Illustrative Example}{Illustrative Example}

\begin{document}

\begin{frontmatter}



\title{Radial Epiderivative Based Line Search Methods in Nonconvex and Nonsmooth Box-Constrained Optimization}


\author[inst1,inst2]{Refail Kasimbeyli}   
\author[inst1]{Gulcin Dinc Yalcin}
\author[inst3]{Gazi Bilal Yildiz}
\author[inst3]{Erdener Ozcetin}

\affiliation[inst1]{Department of Industrial Engineering,
                Eskisehir Technical University, Iki Eylul Campus, Eskisehir 26555, Turkiye}
            
\affiliation[inst2]{UNEC Mathematical Modeling and Optimization Research Center,\\
			   Azerbaijan State University of Economics (UNEC ), Istiqlaliyyet Str. 6, Baku 1001, Azerbaijan.}

\affiliation[inst3]{Department of Industrial Engineering,
              Hitit University, Corum 19030, Turkiye}

\begin{abstract}
In this paper, we develop a novel radial epiderivative-based line search methods for solving nonsmooth and nonconvex box-constrained optimization problems. The rationale for employing the concept of radial epiderivatives is that they provide necessary and sufficient conditions for both identifying global descent directions and achieving global minimum of nonconvex and nondifferentiable functions. These properties of radial epiderivatives are combined with line search methods to develop iterative solution algorithms. The proposed methods generate search directions at each iteration where global descent directions and stopping criteria are performed by using the abilities of the radial epiderivatives. We use two line search methods, that is cyclic coordinate and particle swarm optimization techniques to generate search directions, selecting only those that exhibit descent, as determined by using approximately computed radial epiderivatives at the current point. As a particular case, these methods are applied for minimizing concave functions. In the paper, two convergence theorems are proved. One of them deals with the general line search method and covers only the set of directions generated by the method. The second convergence theorem deals with minimizing concave functions which deals not only with the generated set of directions but covers the whole set of feasible solutions. The performance of the proposed method is evaluated by using well-known benchmark problems from the literature. The results demonstrate the advantages of the proposed approach in generating optimal or near-optimal solutions.

\end{abstract}

\begin{graphicalabstract}
\includegraphics{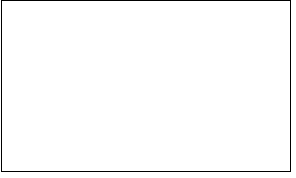}
\end{graphicalabstract}

\begin{highlights}
\item Radial epiderivative based line search methods for solving box-constrained problems.
\item Determining global descent directions for nonconvex and nondifferentiable functions by using radial epiderivatives.
\item Global minimization method for concave functions.
\item Convergence theorems for line search methods to minimize nonconvex and nondifferentiable functions.
\end{highlights}

\begin{keyword}
Radial epiderivative \sep line search methods \sep nonconvex optimization \sep nonsmooth optimization \sep box-constrained optimization \sep global descent direction \sep global solution method 
\PACS 0000 \sep 1111
\MSC 0000 \sep 1111
\end{keyword}

\end{frontmatter}


\section{Introduction}
\label{sec:sample1}
Numerous disciplines, including engineering, finance, and physics, encounter continuous, discrete, constrained or unconstrained, differentiable or nondifferentiable, convex or nonconvex optimization problems. Maximization of profit, reduction of costs, and improvement of manufacturing processes are some examples of issues that can be handled by optimization methods.  

In this paper, we focus on solution methods for the following box-constrained optimization problem:
\begin{equation}\label{problem}
\min_{x \in \mathbb{X}} f(x),
\end{equation}
where $X= \{x=(x_1, x_2, \ldots, x_n): a_i\leq x_i \leq  b_i, i=1, \ldots, n\},$ is a nonempty compact set of feasible solutions, $f : \mathbb{X} \rightarrow \mathbb{R} \cup \{ + \infty \}.$

Algorithms for investigating global and $/$ or local optimal solutions for problems like \eqref{problem}, generally start with an initial solution and need to find an efficient descent direction for the objective function at a current point. Once such a direction is determined, the next point can be generated using a formula like $x_{k+1}=x_k + \alpha_k d_k$ where $k,  \alpha_k$, and $d_k$ denote iteration number, step-size parameter and search direction, respectively.

Generally derivative-based methods use anti-gradient directions as a descent direction (see e.g. \citep{Bazaraa2013}). However, even if the gradient becomes zero (or relatively small) at some point, it may be a local minimum or maximum, depending on whether the function is convex or concave on some neighborhood of this point. When the objective function is nonsmooth, then classical subgradients \citep{Shor} for convex functions or Clarke's subgradients \citep{Clarke} in nonconvex case can be used instead of gradients, to determine search directions; however, again these directions may not be descent in general. For this reason, some methods are developed to compute descent directions by solving (linear or quadratic) subproblems that are formed by using information on (generalized) subgradients.  

In this paper, we develop new radial epiderivative-based line search methods for solving nonsmooth and nonconvex box-constrained optimization problems. In these methods, new search directions are generated by using line search methods where
each direction is investigated whether it is descent or not by using radial epiderivatives. To determine the descent direction $h \in \mathbb{R}^n$ at a given point $\bar{x} \in \mathbb{R}^n$,  we use the radial epiderivative notion $f^r(\overline{x};h)$ introduced by \cite{Kasimbeyli2009}. The radial epiderivative is a concept that extends traditional directional derivative, taking into account changes in a function's values when moving radially in the space of its input variables. The reason for using the radial epiderivative is that, it has become a powerful tool for investigating optimality conditions and descent directions in nonconvex and nonsmooth optimization (see e.g. \citep{DYK2024,KasimbeyliM2009,KasimbeyliM2011}). We use the global optimality condition to formulate stopping criteria and the descent direction condition to determine new search directions in the new algorithms. The proposed method offers several advantages over existing line search methods in the literature. One key benefit is that, by utilizing radial epiderivatives to determine the descent direction, it eliminates the need to solve an additional univariate problem for step length calculation, used in conventional methods. Additionally, we employ an algorithm to approximate the radial epiderivative, which simultaneously computes the step length, further enhancing efficiency.

Many real-life applications are inherently multi-modal, with discontinuities and nonsmooth characteristics, prompting researchers to explore various heuristics and metaheuristics \citep{Wu2015}. The primary advantage of these methods is their ability to solve optimization problems without requiring knowledge of mathematical properties such as convexity and differentiability. Among them, Particle Swarm Optimization (PSO) is one of the most widely used population-based metaheuristics. In PSO, the direction (velocity) of each particle is determined based on its current position, its best-known solution, and the best solution found by the swarm. Due to its strategic approach to direction finding, PSO has been extensively applied to both constrained and unconstrained optimization problems. For instance, Kayhan et al. proposed a hybrid PSO algorithm for continuous optimization \citep{Kayhan2010}, while Fan et al. introduced a hybrid PSO combined with simplex search for unconstrained optimization \citep{Fan2007}. Additionally, \cite{Liang2006} addressed multi-modal functions using a PSO with a comprehensive learning mechanism. More recently, Kwakye et al. enhanced PSO with novel mechanisms for global optimization and feature selection \citep{Kwakye2024}.

In addition to these methods, search directions can also be determined using coordinate axes, as in the Coordinate Cyclic (CC) method (see, e.g., \citep{Bazaraa2013,Wright2015}). CC methods have been extensively studied in the literature for various purposes, including their non-asymptotic convergence behavior \citep{Saha2013}, iteration complexity \citep{Yun2014}, and performance in terms of asymptotic worst-case convergence \citep{Gurbuz2017}. Further research has explored their worst-case complexity for minimizing convex quadratic functions \citep{Sun2021}, local linear convergence in composite nonsmooth optimization problems \citep{Klopfenstein2024}, and extensions to online scenarios where the objective function changes after a finite number of iterations \citep{Lin2024}. Additionally, CC methods have been employed as distributed, model-free reinforcement learning algorithms for designing distributed linear quadratic regulators \citep{Jing2024}, among other applications.

Building on this idea, we decided to combine the strengths of PSO and CC in generating diverse and converging search directions with the ability of the radial epiderivative to determine whether a given direction leads to descent. We then use only those directions that guide us to a better solution. Based on this approach, we developed two new algorithms: Radial Epiderivative-Based PSO (RPSO) and Radial Epiderivative-Based CC (RCC). The performance of these algorithms is evaluated using benchmark test problems from the literature.

The rest of the paper is organized as follows. Section~2 gives definition, and some important properties of radial epiderivative. In this section we also formulate and prove the theorem on concavity for the radial epiderivatives of concave functions. This theorem is used in section~3, to prove the convergence theorem of the proposed method for minimizing concave functions. Radial epiderivative based method is presented in section~3. In this section we present an algorithm for computing radial epiderivatives, whose output is used in Algorithm~2 which explains the general Radial Epiderivative Based Method. Algorithms 3 and 4 explain RPSO and RCC methods respectively. Then, we formulate and prove convergence theorems for the presented algorithms. Computational results for test problems, are given in section ~4. In this section we present explanatory computations and comparisons for the two methods given in this paper. Finally, section~5 draws some conclusions from the paper and an outlook for further research.

\section{Radial Epiderivative} \label{Radialepidersection}

In this section, we recall the definition of the radial epiderivative and outline the algorithm used for its approximate computation. We also present new properties of the radial epiderivative, which will play a key role in the analysis and proof of convergence for the methods developed in this work. 

The radial epiderivative $f^r(\overline{x};\cdot)$ of a given function $f : R^n \rightarrow R$ at a point $\overline{x} \in R^n,$ has been defined by \cite{Kasimbeyli2009} as a function, whose epigraph equals the radial cone to the epigraph $epi(f)$ of $f$ at $(\overline{x}, f(\overline{x}))$ (see also \citep{Flores-Bazan2003}). Kasimbeyli and Mammadov investigated nonsmooth variational problems and optimality conditions by using radial epiderivatives in \citep{KasimbeyliM2009,KasimbeyliM2011}. Recently, \cite{DYK2024} proved that the radial epiderivative $f^r(\overline{x}; d)$ at $\overline{x} \in R^n$ in a direction $d \in R^n,$ can equivalently be defined in the following form:

\begin{equation} \label{radepilimdef}
f^{r}(\overline{x};d) = \inf_{t > 0} \liminf_{ u \rightarrow d} \frac{f(\overline{x}+ tu)- f(\overline{x})}{t}
\end{equation}
for all $\overline{x}, d \in \mathbb{R}^n.$ In the case when the radial epiderivative $f^r (\overline{x}; h)$ exists and finite for every $d,$ we say that $f$ is radially epidifferentiable at $\overline{x}.$

The main advantages of methods developed in this paper are those, they choose descent directions among all the directions considered at every iteration and use the stopping criteria which guarantees global optimality. Both of these benefits are due to the relevan properties of the radial epiderivative notion. The necessary and sufficient condition for the global minimum, given below, was proved by Kasimbeyli in \cite[Theorem 3.6]{Kasimbeyli2009} (see also \cite[Corollary 3]{DYK2024}). 


\begin{theorem} \label{radepiderglobalmin}
	Let $f:R^n \rightarrow \mathbb{R} \cup \{+\infty\}$  be a proper function, radially epidifferentiable at $\overline{x} \in R^n.$ Then $f$ attains its global minimum at  $\overline{x}$ if and only if   $f^r(\overline{x}; \cdot)$ attains its minimum at $d=0.$  
\end{theorem}	


We will say that $d \in R^n$ is a descent direction for function $f:R^n \rightarrow \mathbb{R} \cup \{+\infty\}$ at $\overline{x} \in R^n,$ if there exists a pozitive number $t$ such that $f(\overline{x}+ t d)<f(\overline{x}).$
The following theorem on the necessary and sufficient condition for a descent direction was proved by Dinc Yalcin and Kasimbeyli in  \cite[Theorem 11]{DYK2024}.


\begin{theorem} \label{radepiderdescent}
	Let $f:R^n \rightarrow \mathbb{R} \cup \{+\infty\}$  be a proper function, radially epidifferentiable at $\overline{x} \in R^n.$ Then the vector $d \in R^n$ is a descent direction for $f$ at $\overline{x}$ if and only if   $f^r(\overline{x}; d) < 0.$  
\end{theorem}


In this paper, in addition to the general test problems, we also address minimization problems with concave objective functions. These problems are inherently complex and nonconvex. We apply the proposed methods to concave minimization problems and evaluate their performance on various test cases. In the following section, we present and prove the convergence theorem for the new methods, as well as for the method applied to problems with concave objective functions. Here, we present the following theorem that highlights a key property of the radial epiderivative for concave functions: we demonstrate that the radial epiderivative of a concave function is itself concave.


\begin{theorem}\label{radepiderconcave} Assume that $f:\mathbb{R}^n\rightarrow \mathbb{R}\cup \{+\infty\}$ is a proper  concave function. Then its radial epiderivative $f^r(\overline{x}; \cdot)$ is a concave function of a direction vector, for every $\overline{x} \in R^n.$
\end{theorem}	 
\begin{proof} 
	Let $D=\{d_1, \ldots, d_m \} \subset R^n$ be a given set of directions, $\overline{x} \in R^n$ and let the nonnegative numbers $\{ \alpha_1, \ldots, \alpha_m\} $ satisfy $ \alpha_i \geq 0, i=1, \ldots, m, \sum_{i=1}^{m} \alpha_i =1.$ Then, since $f$ is concave on $R^n,$ it is continuous and hence its radial epiderivative can be written as follows: 
	\begin{eqnarray*}
		f^{r}(\overline{x};\sum_{i=1}^{m} \alpha_i d_i) &=& \inf_{t > 0} \frac{f(\overline{x}+ t\sum_{i=1}^{m} \alpha_i d_i)- f(\overline{x})}{t} = \inf_{t > 0} \frac{f(\sum_{i=1}^{m} \alpha_i (\overline{x}+ td_i))- f(\overline{x})}{t} \\
		&\geq & \inf_{t > 0} \frac{\sum_{i=1}^{m} \alpha_i[ f(\overline{x}+ td_i)- f(\overline{x})]}{t} \geq \sum_{i=1}^{m} \alpha_i \inf_{t > 0} \frac{[ f(\overline{x}+ td_i)- f(\overline{x})]}{t} \\
		&=& \sum_{i=1}^{m} \alpha_i f^{r}(\overline{x}; d_i),
	\end{eqnarray*}
	which proves the assertion of the theorem.	
	
	$\Box$
\end{proof} 



\section{Radial Epiderivative Based Line Search Methods}\label{radialepidermethod}

In this section, we present radial epiderivative-based line search methods for minimizing a given function, without requiring convexity or differentiability conditions. These methods utilize the optimality and descent direction conditions, which are formulated using radial epiderivatives, as discussed in the previous section. To begin, we recall the algorithm for approximately computing radial epiderivatives (see Algorithm \ref{alg: RadialAlg} below), proposed by Dinc Yalcin and Kasimbeyli in \cite[Algorithm 1]{DYK2024}. We use this algorithm to approximate the value of the radial epiderivative, and the resulting output is then applied in the general radial epiderivative-based minimization algorithm. This helps determine whether the generated direction at the current iteration is a descent direction or if the current point is optimal. Finally, this general algorithm is further enhanced by incorporating particle swarm optimization and cyclic coordinate descent methods for generating direction vectors.

Let $\overline{x} \in X$ be a given point and let $h \in R^n$ be a direction vector. Let $f: R^n \rightarrow R$ be a continuous function, lower Lipschitz at $\overline{x}.$ 

\begin{algorithm}
	\caption{Approximate computing of the radial epiderivative $f^r(\overline{x};h)$ of function $f$ at $\overline{x}$ in direction $h.$}
	\label{alg: RadialAlg}
	\begin{algorithmic} [1] 
		\State Let $\overline{x} \in X;$ $t_0>0$ be an initial value for $t$ (which is a sufficiently small positive number); and $ h \in R^n $ be a direction vector, chosen such that $\overline{x} +t_0h \in X$ for $t;$  $\beta >0 $ be a stepsize for $t.$
		\State $y_0 =\frac{f(\overline{x}+t_0h)-f(\overline{x})}{t_0}.$
		\State $\tilde{x}(\overline{x},h) \gets \overline{x}$
		\State $ k \gets 0.$
		\While{$\tilde{x}(\overline{x},h) \in X$ } 
		\State $t_{k+1}= t_0+(k+1)*\beta,$ 
		\State $\tilde{y}_{k+1} =\frac{f(\overline{x}+t_{k+1}h)-f(\overline{x})}{t_{k+1}},$
		\If{$\tilde{y}_{k+1}<y_k$}
		\State $y_{k+1} \gets \tilde{y}_{k+1}.$
		\State $\tilde{x}(\overline{x},h) \gets \overline{x}+t_{k+1}h$
		\Else
		\State $y_{k+1} \gets y_k.$
		\EndIf
		\State $k \gets k+1.$ 
		\EndWhile
		\State $f^r(\overline{x};h)=y_k.$
		\State \Return $\tilde{x}(\overline{x},h).$
	\end{algorithmic}
\end{algorithm}

Algorithm \ref{alg: RadialAlg} explores the function's domain on a grid that is determined by the values of $\overline{x}, h, t_0, \overline{t}$, and $\beta$. During each iteration, the current values for $\tilde{y}_{k+1}$ and $\tilde{x}$ are calculated, and if the smaller ratio is identified, the current approximate value $y_{k+1}$ for the radial epiderivative and the corresponding end point  $\tilde{x}(\overline{x},h) = \overline{x}+t_{k+1}h$  of the ray starting at $\overline{x}$ in the direction of $h,$ are updated. Otherwise, they remain unchanged. The generated point $\tilde{x}(\bar{x},h)$  at the end of the algorithm, is then assigned to the point $x_{k+1}$ in the Algorithm \ref{alg: RadialBasedAlg}.

Dinc Yalcin and Kasimbeyli proved that (see \cite[Theorem 11]{DYK2024}), Algorithm \ref{alg: RadialAlg} computes a finite value for the radial epiderivative $f^r(\overline{x};h)$ after a finite number of iterations, provided that the objective function $f$ is lower Lipschitz at $\overline{x}.$ For the sake of completeness, and because we have made minor changes to both the algorithm and the theorem, we present the theorem along with its proof.


\begin{theorem} \label{algisfinite}
	Assume that $f:\mathbb{R}^n\rightarrow \mathbb{R}\cup \{+\infty\}$ is a proper function, which is lower Lipschitz at $\overline{x} \in X. $ Assume also that the sequences $\{t_k\}$ and $ \{y_k\}$ are generated by Algorithm \ref{alg: RadialAlg} for $k=1,2, \ldots .$ Then, there exists a positive number $N$ such that $y_k = y_N$ for all $k>N.$
\end{theorem}
\begin{proof} 
	Let $f:X \rightarrow \mathbb{R}$ be lower Lipschitz at $\overline{x}$ with Lipschitz constant $L>0.$ Suppose that sequences $\{t_k\}$ and $ \{y_k\}$ are generated by Algorithm \ref{alg: RadialAlg} for $k=1,2, \ldots .$  Assume to the contrary that Algorithm \ref{alg: RadialAlg} generates strongly decreasing sequence of numbers $y_k$ with $y_k >y_{k+1}$ for all $k=1,2, \ldots$ such that this sequence is not bounded from below. Let the number $M$ be chosen such that $M > L.$ Then by the assumption there exists a number $k$ such that $y_k < -M.$ This leads to
	$$
	-M > y_k = \frac{f(\overline{x}+t_kh)-f(\overline{x})}{t_k} \geq \frac{-t_k L \|h\| }{t_k} =-L,
	$$
	which is a contradiction.
	$\Box$
	
\end{proof}


Now we introduce a general global minimization algorithm, based on radial epiderivatives (see Algorithm \ref{alg: RadialBasedAlg}).

\begin{algorithm}
	\caption{Radial Epiderivative based Method}
	\label{alg: RadialBasedAlg}
	\begin{algorithmic} [1] 
		\State Define $x_0, iteration limit, \varepsilon>0$.
        \State Define the direction set $H=\{h_i \in \mathbb{R}^n: \|h_i\|=1 \quad \forall i=\overline{1,m}\}$
		\State Set $distance \gets \infty$
		\State $k \gets 0$
		\While{$distance < \varepsilon$ (or $k \leq iteration limit$ )}
        \State $H_k \gets H$
  		\State Choose a direction $h_k \in H_k.$
		\State Compute the radial epiderivative $f^r(x_k;h^k)$ by using Algorithm \ref{alg: RadialAlg}.
        \While {$f^r(x_k;h_k)>0$ or $H_k \neq \emptyset$}
        \State $H_k \gets H_k \setminus h_k$
        \State Choose a direction $h_k \in H_k.$
        \State Compute the radial epiderivative $f^r(x_k;h_k)$ by using Algorithm \ref{alg: RadialAlg}.
        \EndWhile
        \If{$f^r(x_k;h_k)<0$}
        \State Use the output of  Algorithm \ref{alg: RadialAlg} to define new point $\tilde{x}(x_k,h_k)$ in the direction $h_k$.
        \State $x_{k+1} \gets \tilde{x}(x_k,h_k)$
        \Else
        \State $x_{k+1} \gets x_k$
        \EndIf
		\State $distance \gets  \|x_{k+1}-x_k\|$
		\State $k \gets k+1.$ 
		\EndWhile
		\State \Return $x_k.$
	\end{algorithmic}
\end{algorithm}

In Algorithm \ref{alg: RadialBasedAlg}, the initial steps involve the selection of a starting point $x_k$ and the specific direction $h_k.$ Then the radial epiderivative $f^r(x_k;h_k)$ of $f$ at $x_k$ in the direction $h_k$ is computed by using Algorithm \ref{alg: RadialAlg}. If the value of the radial epiderivative $f^r(x_k;h_k)$ is nonnegative, then an alternative direction is selected. These calculations are repeated until either a direction with a negative radial epiderivative value is identified, or all possible directions have been considered. A new point is generated along the direction with the negative radial epiderivative value, or it is remained unchanged if the radial epiderivative values are positive in all possible directions. Subsequently, if the distance between the newly generated point and the current point becomes less than a predetermined accuracy value $\varepsilon$, this indicates that there is no significant decrease in the objective function or that there is no descent direction. In this situation, the algorithm terminates. Otherwise, the next iteration is performed. 

The direction set used in Algorithm \ref{alg: RadialBasedAlg} can be generated using various strategies. In this paper, we employ two methods for generating the set of possible directions. First, we use the Particle Swarm Optimization (PSO) method, a widely used population-based metaheuristic. Next, we implement the cyclic coordinate descent method to generate the direction vectors.\\


~\\

\subsection{Radial epiderivative based method with particle swarm optimization}

In this subsection, we introduce the radial epiderivative based method where the direction set is generated by using Particle Swarm Optimization (PSO) algorithm. PSO has been used extensively in domains including engineering, machine learning, and operational research because of its ease of use, effectiveness, and capacity to manage nonlinear and multidimensional issues. This is a population-based stochastic algorithm inspired by a nature that solves challenging optimization problems by imitating the social behavior of bird flocking. It was first presented by \cite{kennedy1995particle}. Our goal is to combine the (PSO) algorithm to generate a set of feasible directions with the capabilities of radial epiderivatives to select only descent directions. Additionally, we aim to apply the necessary and sufficient conditions for global optimality using radial epiderivatives.

PSO employs a set of particles ($i$) to explore the search space, with each particle representing a potential solution ($x_k^i$). These particles gradually converge towards a (possibly) better feasible solution by modifying their positions in response to their own experiences ($x_{pbest}^i$), those of their neighbours ($x_{gbest}$), and the velocity ($v$) of the particles combined with random variables $r_1$ and $r_2$ which fall within the range of $[0,1]$ and with the parameters $w, c_1$ and $c_2$ of the algorithm. In this paper, the idea of PSO is applied to generate the direction set by using the velocity (direction) of each particle at iteration $k+1$ as follows:
\begin{equation}\label{velocity}
v_{k+1}^i= w v_k^i+c_1r_1(x_{pbest}^i-x_k^i)+c_2r_2(x_{gbest}-x_k^i).
\end{equation}
The determination of the suitability of the velocity as a descent direction is established by employing the radial epiderivative value at the current point of the objective function in each iteration. In the case of a descent direction, the position of each particle undergoes an update; in the absence of such a descent direction, the particle's position remains unchanged. The method that entails the utilization of radial epiderivative in conjunction with Particle Swarm Optimization (PSO) is presented in Algorithm \ref{alg: RadialBasedAlgPSO}.

\begin{algorithm}
	\caption{Radial Epiderivative based Method with Particle Swarm Optimization (RPSO)}
	\label{alg: RadialBasedAlgPSO}
	\begin{algorithmic} [1] 
		\State Define the number of particles $P,$ the position of each  particle $i$ as $x_0^i$.
        \State $iteration limit, \varepsilon,$ and $\overline{count}.$
        \State Define parameters $t_0, \underline{t}, \alpha \in (0,1), \beta$ for Algorithm 1.
        \State Set $t_0^i \gets t_0, \beta^i \gets \beta \quad \forall i \in P$.
        \State Define the minimum value $\underline{t}$ for $t_0.$
		\State $x_{pbest}^i \gets x_0^i \quad \forall i \in P, x_{gbest} \gets \argmin_i\{f(x_0^i)\}.$
		\State $k \gets 0.$
		\While{$(k \leq iteration limit)$ or $(count < \overline{count})$} 
		\For{$i \in P$}
		\State Copmute the velocity of the particle $v_k^i$ by using \eqref{velocity}.
        \If{$x_{gbest}$ and $x_{pbest} ^i$ is not changed at the previous iteration $k-1$}
        \If {$t_0^i>\underline{t}$}
        \State $t_0^i \gets t_0^i*\alpha$
        \State $\beta^i \gets \beta^i* \alpha$
        \Else
        \State Generate the position of of  particle $i$ as $x_{k+1}^i$ randomly
        \State $t_0^i \gets t_0$
        \State $\beta^i \gets \beta$
        \EndIf
        \Else
        \If{$x_{pbest} ^i$ is not changed at the previous iteration $k-1$}
        \If {$t_0^i>\underline{t}$}
        \State $t_0^i \gets t_0^i*\alpha$
        \State $\beta^i \gets \beta^i* \alpha$
        \EndIf
        \EndIf
        \EndIf
		\State Compute radial epiderivative $f^r(x_k^i;v_k^i)$ by using Algorithm \ref{alg: RadialAlg}.
        \State $x_{k+1}^i \gets \Tilde{x}$
        \If {$f(x_{k+1}^i) < f(x_{pbest}^i)$}
        \State $x_{pbest} ^i \gets x_{k+1}^i $
        \EndIf
		\EndFor
        \If {$\min_i\{f(x_{k+1}^i)\}  < f(x_{gbest})$}
        \State $x_{gbest} \gets \argmin_i\{f(x_{k+1}^i)\} $
        \EndIf
        \If{$x_{gbest}$ is not changed}
		\State $count \gets count +1$
        \Else
        \State $count \gets 0$
        \EndIf
		\State $k \gets k+1.$ 
		\EndWhile
		\State \Return $x_{gbest}$ 
	\end{algorithmic}
\end{algorithm}

In Algorithm \ref{alg: RadialBasedAlgPSO}, due to the presence of multiple particles, nonupdated global best values are counted. Subsequently, in the case that both the global best solution and the particle's best solution have not been updated in the previous iteration and the ratio of $t_0^i$ is relatively small for the particle, the position of the particle is randomly generated once again. This approach is justified as it implies that the velocity remains constant and a grid search is conducted within a sufficiently small ratio. Consequently, it becomes unnecessary to conduct a search in the same area any further. If the global best value remains the same after a predetermined number of consecutive iterations $\overline{count}$ or the maximum iteration number is reached, then the algorithm terminates. 

\subsection{Radial epiderivative based method with cyclic coordinate}

The direction in Algorithm \ref{alg: RadialBasedAlg} can be figured out by employing the coordinate directions defined as basis vectors  $d_i = (0,0,\ldots,1,...,0)^T, \mbox{ for } i \in \{1,\ldots,n\}$ where $1$ stands ate the $i^th$ position, or $-d_i, \mbox{ for } i \in \{n+1,...,2n\}.$  In the case if the radial epiderivative associated with the direction $d_i$ (or $-d_i$) has a negative value, this indicates that the direction represents a descent, leading to an update of the current point; conversely, if the value is non-negative, the current point remains unchanged. We elucidate the radial epiderivative based technique with cyclic coordinate in Algorithm \ref{alg: RadialBasedAlgCC}. 

\begin{algorithm}  
	\caption{Radial Epiderivative based Method with Cyclic Coordinate (RCC)}
	\label{alg: RadialBasedAlgCC}
	\begin{algorithmic} [1] 
		\State Define $x_0, iteration limit$ and $\varepsilon>0$
		\State Set $distance \gets \infty$
		\State Let $d_i = (0,0,\ldots, 1,...,0)^T, \mbox{ or }  d_i = (0,0,\ldots, -1,...,0)^T \mbox{ for } i \in \{1,...,2n\}$
		\State $k \gets 0$
		\While{$(k \leq iteration limit)$ or $(count < \overline{count})$} 
        \State $f_{best} \gets f(x_k)$
        \State $x_{best} \gets x_k$
		\For {$j \leq 2n$}
		\State Compute radial epiderivative $f^r(x_k;d_j)$ by using Algorithm \ref{alg: RadialAlg}.
		\State $y_{j} \gets \overline{x}$
        \If{$f(y_j)<f_{best}$}
        \State $f_{best} \gets f(y_j)$
        \State $x_{best} \gets y_j$
        \EndIf
        \EndFor
        \State $x_{k+1} \gets x_{best}$
		\State $distance \gets  \|x_{k+1}-x_k\|$
        \If{$distance < \varepsilon$ }
        \State $t_0 \gets t_0*\alpha$
        \State $\beta \gets \beta* \alpha$
        \State $count \gets count +1$
        \EndIf
		\State $k \gets k+1.$ 
		\EndWhile
		\State \Return $x_k$  
	\end{algorithmic}
\end{algorithm}

As delineated in Algorithm \ref{alg: RadialBasedAlgCC}, the implementation of Algorithm \ref{alg: RadialBasedAlg} is executed by utilizing the directions of basic vectors $d_i$ at each iteration. The evaluation of the radial epiderivative of the function at the relevant point along a specified direction of basis vectors ascertains whether that direction constitutes a descent direction. Upon the comprehensive examination of all orientations of basis vectors, the metric of separation between the current point and the subsequently generated point is computed.  If this metric is inferior to a pre-established threshold, $\varepsilon$, the grid interval utilized for the computation of the radial epiderivative is subsequently diminished.
Consequently, the algorithm advances in a comparable fashion. If the computed distance is found to be less than the threshold value of $\varepsilon$ after a specified number of iterations denoted as $\overline{count},$ or if the iteration limit is attained, the algorithm is concluded.

Next section provides some computational results demonstrating numerical performance of the proposed method. Now we prove theorems characterizing convergence properties of the sequences of feasible solutions generated by the above algorithms.

\begin{theorem}(Convergence Theorem)\label{convergmain} Consider the problem \eqref{problem}.
	Let $f:\mathbb{R}^n\rightarrow \mathbb{R}\cup \{+\infty\}$ be a proper  lower semi-continuous function and let the sequence $\{ x_k \}$ be generated by Algorithm \ref{alg: RadialBasedAlgCC} in the form $x_{k+1} = x_k + t_k d_k, k = 1, 2, \ldots$ where the direction $d_k$ is selected from a set of vectors $D=\{d_1, \ldots, d_m \},$ where $m \geq 1$ is an arbitrary natural number, according to the descent direction choosing rule by using the radial epiderivative test $f^r(x_k;d_k) <0.$ If there are more than one descent direction at a point $x_k,$ the direction that provides the smallest value to the objective function is selected. Asssume that Algorithm \ref{alg: RadialBasedAlgCC} generates an infinite sequence $\{ (x_k,d_k) \}, k = 1, 2, \ldots , .$ Then, every cluster point $(\overline{x},\overline{d}) \in X \times D$ of the sequence $\{ (x_k,d_k) \}$ satisfies $ 0 = f^r(\overline{x};\overline{d}) \leq f^r(\overline{x};d)$ for every $d \in D.$
\end{theorem}	 
\begin{proof} Assume that the sequence $\{ x_k \}$ is generated in the form $x_{k+1} = x_k + t_k d_k$ with $f^r(x_k;d_k) < 0$ and hence $f(x_{k+1}) < f(x_k).$
Since both the sets  $X$ and $D$ are compact, at least one cluster point of the sequence  $\{ (x_k,d_k) \}$ must exist. Without loss of generality, assume that $ lim_{k\rightarrow \infty} x_k =\overline{x} $ and $lim_{k\rightarrow \infty} d_k =\overline{d} \in D. $ 
Assume by contradiction that $f^r(\overline{x};\overline{d}) < 0.$ This means that $\overline{d}$ is a descent direction at $\overline{x}$ and therefore there exists a number $\overline{t}>0$ such that $\delta\colon= f(\overline{x}) - f(\overline{x} + \overline{t} \overline{d}) >0$ and  $\overline{x} + \overline{t} \overline{d} \in int X.$ Then since $(x_k + \overline{t} d_k) \rightarrow (\overline{x} + \overline{t} \overline{d}), $ for $k$ sufficiently large we have:
\begin{equation*}
f(x_k + \overline{t} d_k) \leq f(\overline{x} + \overline{t} \overline{d}) + \frac{\delta}{2} = f(\overline{x}) - \delta + \frac{\delta}{2} = f(\overline{x}) - \frac{\delta}{2}.
\end{equation*}
However,
\begin{equation*}
f(\overline{x}) <    f(x_k + t_k d_k) \leq f(x_k + \overline{t} d_k) \leq  f(\overline{x}) - \frac{\delta}{2},
\end{equation*}
which is a contradiction. Thus $f^r(\overline{x};\overline{d}) = 0.$

Now we show that there does not exist a direction $d \in D$ with $f^r(\overline{x};d) < 0.$ Again assume by contradiction that such a direction $\tilde{d}$ does exist, and hence there exists a number $\tilde{t}>0$ such that $f(\overline{x} + \tilde{t} \tilde{d}) < f(\overline{x}).$ Let $f(\overline{x}) - f(\overline{x} + \tilde{t} \tilde{d}) = \delta >0 $ and $\overline{x} + \tilde{t} \tilde{d} \in int X.$  Then since $ x_k + \tilde{t} \tilde{d} \rightarrow \overline{x}  + \tilde{t} \tilde{d} $ and since $ f(x_k + \tilde{t} \tilde{d}) \geq f(x_k +t_kd_k) \geq f(\overline{x}) > f(\overline{x} + \tilde{t} \tilde{d}), $ for $k$ sufficiently large we have:
\begin{equation*}
f(\overline{x}) \leq f(x_k +t_kd_k) \leq f(x_k + \tilde{t} \tilde{d}) \leq f(\overline{x} +  \tilde{t} \tilde{d}) + \frac{\delta}{2} = f(\overline{x}) - \delta + \frac{\delta}{2} = f(\overline{x}) - \frac{\delta}{2}.
\end{equation*}
Which is a contradiction and hence the proof is completed.
$\Box$
\end{proof} 


\begin{corollary}\label{alg4genglobalmin} 
	Theorem \ref{convergmain} proves that if Algorithm \ref{alg: RadialBasedAlgCC} generates an infinite sequence $\{ (x_k,d_k) \},$ for $ k = 1, 2, \ldots ,$ then, every cluster point $(\overline{x},\overline{d})$ of the sequence $\{ (x_k,d_k) \}$ provides a smallest value to $f$ among the all directions considered by the algorithm.
\end{corollary}

Theorem \ref{convergmain} and Corollary \ref{alg4genglobalmin} show that Algorithm \ref{alg: RadialBasedAlgCC} the performance of the method can be improved by extending the set of feasible directions considered by the algorithm. From this point of view, we can expect a better performance from the Radial Epiderivative based Method with Particle Swarm Optimization (RPSO). On the other hand, the following theorem shows that the sequence of feasible solutions generated by the method presented in this paper,  converges to a global minimum of the problem \eqref{problem}, if the objective function $f$ is a proper concave function. It is known that minimizing a concave function is a nonconvex problem, often considered one of the most challenging problems in the literature.


\begin{theorem}\label{convergenceforconcavefunc} 
	Consider the problem \eqref{problem} and assume that the objective function  $f:\mathbb{R}^n\rightarrow \mathbb{R}\cup \{+\infty\}$ is a proper concave function. Let the sequence $\{ x_k \}$ be generated by Algorithm \ref{alg: RadialBasedAlgCC} in the form $x_{k+1} = x_k + t_k d_k, k = 1, 2, \ldots$ where the direction $d_k$ is selected from a set of basis vectors $D=\{d_1, \ldots, d_n \}$ or from  $-D=\{-d_1, \ldots, -d_n \},$ according to the descent direction choosing rule by using the radial epiderivative test $f^r(x_k;d_k) < 0.$ If there are more than one descent direction among the all directions from $D \cup -D$  at a point $x_k,$ the direction that provides the smallest value to the objective function is selected. Asssume that Algorithm \ref{alg: RadialBasedAlgCC} generates an infinite sequence $\{ (x_k,d_k) \}, k = 1, 2, \ldots , .$  Then, every cluster point $(\overline{x},\overline{d})$ of the sequence $\{ (x_k,d_k) \}$ satisfies $f^r(\overline{x};\overline{d}) = 0$ where $\overline{x}$ is a global minimum of $f$ over $X.$
\end{theorem}	 
\begin{proof} 
	The first part of the proof is the same as in the proof of the proof of Theorem \ref{convergmain}. By using the same way with the set of basis vectors $D=\{d_1, \ldots, d_n \}$ and  $-D=\{-d_1, \ldots, -d_n \},$ we obtain that every cluster point $(\overline{x},\overline{d})$ of the sequence $\{ (x_k,d_k) \}$ satisfies $f^r(\overline{x};\overline{d}) = 0$ and that  $f^r(\overline{x};d_i) \geq 0$ for all $d_i \in D\cup -D.$
	
 Now we show that there does not exist a feasible direction $\tilde{d} \in R^n$ (that is a direction $\tilde{d}$ with $\overline{x} + t \tilde{d} \in int(X)$ for some $t>0$ ) with $f^r(\overline{x};\tilde{d}) < 0.$ Assume by contradiction that such a direction $\tilde{d} \in R^n$ does exist. Then, since the set of vectors $D$ is a basis for $R^n,$ there exists a set of numbers $\lambda_1,\lambda_2, \ldots , \lambda_n $ such that $\tilde{d} = \sum_{i=1}^{n} \lambda_i d_i.$ Without loss of generality we can assume that all the numbers $\lambda_i$ are nonnegative, because if for example $\lambda_j$ is negative, then the corresponding vector $d_j$ in the representation of  $\tilde{d}$ can be changed to $-d_j$ and $\lambda_j$ to $-\lambda_j>0.$ Clearly the new set of vectors $\{d_1, \ldots, d_{j-1}, -d_j, d_{j+1} , \ldots, d_n\}$ is also a basis for $R^n.$ Then we have:
  \begin{equation}
 	0> f^r(\overline{x} ;\tilde{d}) = f^r(\overline{x} ; \sum_{i=1}^{n} \lambda_i d_i).
  \end{equation}
  By multiplying both sides of the above inequality by a pozitive number $\frac{1}{\sum_{i=1}^{n} \lambda_i}$ and using the pozitive homogeneity and concavity of the radial epiderivative, we obtain:
\begin{eqnarray*}
  	0 &>& \frac{1}{\sum_{i=1}^{n} \lambda_i}f^r(\overline{x} ;\tilde{d}) = \frac{1}{\sum_{i=1}^{n} \lambda_i} f^r(\overline{x} ; \sum_{i=1}^{n} \lambda_i d_i) =  f^r(\overline{x} ; \sum_{i=1}^{n} \frac{\lambda_i}{\sum_{i=1}^{n} \lambda_i} d_i) \\
  	& \geq & \sum_{i=1}^{n} \frac{\lambda_i}{\sum_{i=1}^{n} \lambda_i} f^r(\overline{x} ;  d_i) \geq 0,
\end{eqnarray*}
which is a contradiction and hence the proof is completed. 
	$\Box$
\end{proof}



\section{Computational Results}
In this section, we perform computations applying the three algorithms explained in the previous section. In addition, we discuss the strength of the radial epiderivative to find better directions in detail. For computational experiments, we use well-known $29$ different test problems from \cite{AlRoomi2015}. These test problems and their optimal solutions are given in Table \ref{tab:t1} and Table \ref{tab:optimal_solutions}, respectively. These problems are categorized as continuous, differentiable, and convex or nonconvex. In order to investigate the effect of the proposed methods on minimizing concave functions, we have generated several test problems by taking some convex problems from the literature and multiplying the objective functions of these problems by $(-1).$ The reason is that we were unable to find a compelling set of test instances on minimizing concave functions.

The parameters are defined as follows.
\begin{enumerate}
	\item Starting points $x_0$ are selected randomly in PSO and RPSO, and it is defined as the middle point of the lower and upper limits in CC and RCC;
	\item $t_0, \alpha,$ and $\beta$  are all put equal to $0.1$ for the RPSO and the RCC;
	\item $\overline{count}$ is put equal to $3$ for all algorithms.
	\item $f$ denotes the value of the objective function computed by one of the algorithms PSO, RPSO, CC, RCC and $f^*$ denotes the optimal (or best known) value of the objective function. 
\end{enumerate}
 
 All experiments in this section were performed on a computer with an Intel i5-2.7 Ghz. processor with 16 GB RAM. The implementation of algorithms is coded in Python  3.8.5.

Performance measurement is required to demonstrate the performance of the algorithms in comparison to the optimal solutions. In this study, we calculate the gap in the following form (see e.g. \citep{DYK2021,DY2022}):
	
 \begin{equation}\label{EvalCri}
	Gap = \frac{f - f^*}{1 + |f^*|}. 
	\end{equation}

\begin{sidewaystable}
\caption{Test Problems}
\label{tab:t1}
\resizebox{\textwidth}{!}{
\begin{tabular}{lp{1.8cm}p{2cm}p{1.5cm}p{3cm}p{3cm}p{5cm}p{3cm}}
\hline
\textbf{Problem Name} & \textbf{Continuous} & \textbf{Differentiable} & \textbf{Convex} & \textbf{Search Range ($x_1$)} & \textbf{Search Range ($x_2$)} & \textbf{Optimal Solution} & \textbf{Optimum Value} \\ \hline

Ackley1 & + & + & - & [-35, 35] & [-35, 35] & (0, 0) & 0 \\
Alpine1 & + & - & - & [-10, 10] & [-10, 10] & (0, 0) & 0 \\
Brent & + & + & + & [-10, 10] & [-10, 10] & (-10, -10) & $1.3839 \times 10^{-87}$ \\
Brown & + & + & + & [-1, 4] & [-1, 4] & (0, 0) & 0 \\
ChungReynolds & + & + & + & [-100, 100] & [-100, 100] & (0, 0) & 0 \\
Csendes & + & - & - & [-2, 2] & [-2, 2] & (0, 0) & 0 \\
Deb1 & + & + & - & [-1, 1] & [-1, 1] & (-0.1, -0.1) & -1 \\
Deb2 & + & + & - & [0, 1] & [0, 1] & (0.079699, 0.079699) & -1 \\
DixonPrice & + & + & + & [-10, 10] & [-10, 10] & $(1, \frac{1}{\sqrt{2}})$ & 0 \\
DropWave & + & + & - & [-5.12, 5.12] & [-5.12, 5.12] & (0, 0) & -1 \\
EggHolder & + & + & - & [-512, 512] & [-512, 512] & (512, 404.2319) & -959.6407 \\
Exponential & + & + & + & [-1, 1] & [-1, 1] & (0, 0) & 1 \\
Giunta & + & + & - & [-1, 1] & [-1, 1] & (0.46732, 0.46732) & 0.06447 \\
Mishra1 & + & + & - & [0, 1] & [0, 1] & (1, 1) & 2 \\
Mishra2 & + & + & - & [0, 1] & [0, 1] & (1, 1) & 2 \\
Periodic & + & + & - & [-10, 10] & [-10, 10] & (0, 0) & 0.9 \\
PowellSum & + & + & + & [-1, 1] & [-1, 1] & (0, 0) & 0 \\
Qing & + & + & - & [-500, 500] & [-500, 500] & $(1, \sqrt{2})$ & 0 \\
Rastrigin & + & + & - & [-5.12, 5.12] & [-5.12, 5.12] & (0, 0) & 0 \\
Rosenbrock & + & + & - & [-5, 10] & [-5, 10] & (1, 1) & 0 \\
Salomon & + & + & - & [-100, 100] & [-100, 100] & (0, 0) & 0 \\
SchumerSteiglitz & + & + & + & [-100, 100] & [-100, 100] & (0, 0) & 0 \\
Sphere & + & + & + & [-5.12, 5.12] & [-5.12, 5.12] & (0, 0) & 0 \\
Step & - & - & - & [-100, 100] & [-100, 100] & (0, 0) & 0 \\
StepInt & - & - & - & [-5.12, 5.12] & [-5.12, 5.12] & (-5.12, -5.12) & 13 \\
SumSquares & + & + & + & [-10, 10] & [-10, 10] & (0, 0) & 0 \\
Trid & + & + & - & [-8, 8] & [-8, 8] & (2, 2) & -2 \\
Vincent & + & + & - & [0.25, 10] & [0.25, 10] & (7.70628098, 7.70628098) & -2 \\
WWavy & + & + & - & [-$\pi$, $\pi$] & [-$\pi$, $\pi$] & (0, 0) & 0 \\

\hline
\end{tabular}
}
\end{sidewaystable}


~\\

\begin{table}[h!]
	\centering
	\caption{Optimal Solutions of Problems}
	\label{tab:optimal_solutions}
	\renewcommand{\arraystretch}{0.9} 
	\setlength{\tabcolsep}{6pt} 
	\small 
	\begin{tabular}{lcc}
		\toprule
		\textbf{Problem} & \textbf{$x_1$} & \textbf{$x_2$} \\
		\midrule
		\textbf{Ackley1}          & -0.0004 & -0.0008 \\
		\textbf{Alpine1}          & -0.0007 & -0.1002 \\
		\textbf{Brent}            & -9.9994 & -9.9999 \\
		\textbf{Brown}            & 0.0010  & -0.0007 \\
		\textbf{ChungReynolds}    & -0.0012 & 0.0015  \\
		\textbf{Csendes}          & 0.0000  & 0.0000  \\
		\textbf{Deb1}             & 0.5000  & 0.3000  \\
		\textbf{Deb2}             & 0.0797  & 0.0798  \\
		\textbf{DixonPrice}       & 0.9976  & -0.7062 \\
		\textbf{DropWave}         & 0.0000  & 0.0000  \\
		\textbf{EggHolder}        & 512.0000 & 404.2279 \\
		\textbf{Exponential}      & 0.0000  & 0.0000  \\
		\textbf{Giunta}           & 0.0000  & 0.0000  \\
		\textbf{Mishra1}          & 1.0000  & 0.5143  \\
		\textbf{Mishra2}          & 1.0000  & 0.9996  \\
		\textbf{Periodic}         & 0.0002  & 0.0002  \\
		\textbf{PowellSum}        & 0.0000  & 0.0000  \\
		\textbf{Qing}             & 0.9994  & 1.4144  \\
		\textbf{Rastrigin}        & 0.0001  & 0.0014  \\
		\textbf{Rosenbrock}       & 0.9955  & 0.9915  \\
		\textbf{Salomon}          & 0.0005  & 0.0023  \\
		\textbf{SchumerSteiglitz} & -0.0018 & -0.0022 \\
		\textbf{Sphere}           & -0.0036 & 0.0007  \\
		\textbf{Step}             & 0.0326  & -0.8016 \\
		\textbf{StepInt}          & -5.1189 & -5.0002 \\
		\textbf{SumSquares}       & -0.0007 & 0.0018  \\
		\textbf{Trid}             & 1.9974  & 1.9993  \\
		\textbf{Vincent}          & 1.4358  & 1.4359  \\
		\textbf{WWavy}            & 0.0000  & 0.0000  \\
		\bottomrule
	\end{tabular}
\end{table}


~\\

For benchmarking purposes, we evaluated all test problems using not only the proposed methods but also the CC and PSO algorithms. Table \ref{tab:t2} presents the computational results obtained by applying the CC, RCC, PSO, and RPSO methods. The table also includes the gap values associated with each approach. As observed, the RCC method outperforms CC in 24 out of 29 instances, while both methods yield equivalent results in 4 instances. CC surpasses RCC in only one case. A similar trend is evident in the comparison between PSO and RPSO: PSO outperforms RPSO in just one instance, while both methods show identical performance in 12 cases. RPSO demonstrates superior results in the remaining 16 problems. Furthermore, when compared to the optimal solutions, only four results obtained by CC and five by RCC achieve near-optimality with minimal optimality gaps.

\begin{sidewaystable}
\caption{Computational Results of CC, RCC, PSO, and RPSO}
\label{tab:t2}
\centering
\resizebox{\columnwidth}{!}{
\begin{tabular}{llllllllll}
\hline
\multirow{2}{*}{Problems} & \multirow{2}{*}{\textbf{Optimum Value}} & \multicolumn{4}{c}{\textbf{Objective Values}}                                 & \multicolumn{4}{c}{\textbf{Gaps}}                                         \\
                          &                                         & \textbf{CC-Obj}  & \textbf{RCC-Obj}   & \textbf{PSO-Obj} & \textbf{RPSO-Obj}  & \textbf{CC-Gap} & \textbf{RCC-Gap} & \textbf{PSO-Gap} & \textbf{RPSO-Gap} \\
                          \hline
\textbf{Ackley1}          & 0.0000                                  & 19.8794          & \textbf{0.0000}    & 0.0152           & \textbf{0.0025}    & 4.8394          & \textbf{0.7966}  & 0.0152           & \textbf{0.0025}   \\
\textbf{Alpine1}          & 0.0000                                  & 6.4402           & \textbf{0.0000}    & \textbf{0.0001}  & \textbf{0.0001}    & 6.4402          & \textbf{0.0000}  & \textbf{0.0001}  & \textbf{0.0001}   \\
\textbf{Brown}            & 0.0000                                  & 0.0000           & \textbf{0.0000}    & \textbf{0.0000}  & \textbf{0.0000}    & 0.0000          & \textbf{0.0000}  & \textbf{0.0000}  & \textbf{0.0000}   \\
\textbf{ChungReynolds}    & 0.0000                                  & na               & \textbf{0.0000}    & \textbf{0.0000}  & \textbf{0.0000}    & na              & \textbf{0.0000}  & \textbf{0.0000}  & \textbf{0.0000}   \\
\textbf{Csendes}          & 0.0000                                  & 97.3154          & \textbf{0.0000}    & \textbf{0.0000}  & \textbf{0.0000}    & 97.3154         & \textbf{0.0000}  & \textbf{0.0000}  & \textbf{0.0000}   \\
\textbf{Deb1}             & -1.0000                                 & -1.0000          & \textbf{-1.0000}   & -0.9997          & \textbf{-1.0000}   & \textbf{0.0000} & \textbf{0.0000}  & 0.0001           & \textbf{0.0000}   \\
\textbf{Deb2}             & -1.0000                                 & \textbf{-1.0000} & \textbf{-1.0000}   & -0.9998          & \textbf{-1.0000}   & \textbf{0.0000} & \textbf{0.0000}  & 0.0001           & \textbf{0.0000}   \\
\textbf{DixonPrice}       & 0.0000                                  & \textbf{0.0000}  & \textbf{0.0000}    & 0.0001           & \textbf{0.0000}    & \textbf{0.0000} & \textbf{0.0000}  & 0.0001           & \textbf{0.0000}   \\
\textbf{DropWave}         & -1.0000                                 & -0.1077          & \textbf{-1.0000}   & -0.9978          & \textbf{-1.0000}   & 0.4462          & \textbf{0.0000}  & 0.0011           & \textbf{0.0000}   \\
\textbf{EggHolder}        & -959.6407                               & -603.5642        & \textbf{-894.5789} & -114.8907        & \textbf{-959.6405} & 0.3707          & \textbf{0.0677}  & 0.8794           & \textbf{0.0000}   \\
\textbf{Exponential}      & -1.0000                                 & -0.6065          & \textbf{-1.0000}   & \textbf{-1.0000} & \textbf{-1.0000}   & 0.1967          & \textbf{0.0000}  & \textbf{0.0000}  & \textbf{0.0000}   \\
\textbf{Giunta}           & 0.0645                                  & 0.2080           & \textbf{0.0645}    & \textbf{0.0645}  & \textbf{0.0645}    & 0.1348          & \textbf{0.0000}  & \textbf{0.0000}  & \textbf{0.0000}   \\
\textbf{Mishra1}          & 2.0000                                  & \textbf{2.0000}  & \textbf{2.0000}    & \textbf{2.0000}  & \textbf{2.0000}    & \textbf{0.0000} & \textbf{0.0000}  & \textbf{0.0000}  & \textbf{0.0000}   \\
\textbf{Mishra2}          & 2.0000                                  & 2.7539           & \textbf{2.0000}    & \textbf{2.0000}  & 2.0004             & 0.2513          & \textbf{0.0000}  & \textbf{0.0000}  & 0.0001            \\
\textbf{Periodic}         & 0.9000                                  & 1.2960           & \textbf{0.9000}    & \textbf{0.9000}  & \textbf{0.9000}    & 0.2084          & \textbf{0.0000}  & \textbf{0.0000}  & \textbf{0.0000}   \\
\textbf{PowellSum}        & 0.0000                                  & 1.0000           & \textbf{0.0000}    & \textbf{0.0000}  & \textbf{0.0000}    & 1.0000          & \textbf{0.0000}  & \textbf{0.0000}  & \textbf{0.0000}   \\
\textbf{Qing}             & 0.0000                                  & \textbf{0.0000}  & \textbf{0.0000}    & 0.0001           & \textbf{0.0000}    & \textbf{0.0000} & \textbf{0.0000}  & 0.0001           & \textbf{0.0000}   \\
\textbf{Rastrigin}        & 0.0000                                  & 32.9043          & \textbf{0.0000}    & 0.1303           & \textbf{0.0004}    & 32.9043         & \textbf{0.0000}  & 0.1303           & \textbf{0.0004}   \\
\textbf{Rosenbrock}       & 0.0000                                  & \textbf{0.4448}  & 0.6225             & 0.0005           & \textbf{0.0000}    & 0.4448          & 0.6225           & 0.0005           & \textbf{0.0000}   \\
\textbf{Salomon}          & 0.0000                                  & 10.6999          & \textbf{0.0999}    & 0.0107           & \textbf{0.0003}    & 10.6999         & \textbf{0.0999}  & 0.0107           & \textbf{0.0003}   \\
\textbf{SchumerSteiglitz} & 0.0000                                  & na               & \textbf{0.0000}    & 0.0000           & \textbf{0.0000}    & na              & \textbf{0.0000}  & \textbf{0.0000}  & \textbf{0.0000}   \\
\textbf{Sphere}           & 0.0000                                  & 26.2144          & \textbf{0.0000}    & \textbf{0.0000}  & \textbf{0.0000}    & 26.2144         & \textbf{0.0000}  & \textbf{0.0000}  & \textbf{0.0000}   \\
\textbf{Step}             & 0.0000                                  & 101.0000         & \textbf{0.0000}    & \textbf{0.0000}  & \textbf{0.0000}    & 101.0000        & \textbf{0.0000}  & \textbf{0.0000}  & \textbf{0.0000}   \\
\textbf{StepInt}          & 13.0000                                 & 14.0000          & \textbf{13.0000}   & \textbf{13.0000} & \textbf{13.0000}   & 0.0714          & \textbf{0.0000}  & \textbf{0.0000}  & \textbf{0.0000}   \\
\textbf{SumSquares}       & 0.0000                                  & 199.9998         & \textbf{0.0000}    & \textbf{0.0000}  & \textbf{0.0000}    & 199.9998        & \textbf{0.0000}  & \textbf{0.0000}  & \textbf{0.0000}   \\
\textbf{Trid}             & -2.0000                                 & 1.0184           & \textbf{-2.0000}   & \textbf{-2.0000} & \textbf{-2.0000}   & 1.0061          & \textbf{0.0000}  & \textbf{0.0000}  & \textbf{0.0000}   \\
\textbf{Vincent}          & -2.0000                                 & \textbf{-2.0000} & \textbf{-2.0000}   & -1.9998          & \textbf{-2.0000}   & \textbf{0.0000} & \textbf{0.0000}  & 0.0001           & \textbf{0.0000}   \\
\textbf{WWavy}            & 0.0000                                  & 0.4964           & \textbf{0.0000}    & \textbf{0.0000}  & \textbf{0.0000}    & 0.4964          & \textbf{0.0000}  & \textbf{0.0000}  & \textbf{0.0000}   \\
\textbf{Brent}            & 0.0000                                  & \textbf{0.0000}  & \textbf{0.0000}    & \textbf{0.0000}  & \textbf{0.0000}    & \textbf{0.0000} & \textbf{0.0000}  & \textbf{0.0000}  & \textbf{0.0000}  
\\
\hline

\end{tabular}



  %
  
}
\newline
\end{sidewaystable}

\begin{table}
\caption{Computational results based on iteration limits of PSO and RPSO in terms of Gaps }
\label{tab:t3}
\centering
\resizebox{\columnwidth}{!}{
\begin{tabular}{lcccccc}
\hline
\multicolumn{1}{c}{\multirow{2}{*}{\textbf{Problem}}} & \multicolumn{2}{c}{\textbf{iter=100}} & \multicolumn{2}{c}{\textbf{iter=500}} & \multicolumn{2}{c}{\textbf{iter=1000}} \\
\multicolumn{1}{c}{}                                   & \textbf{PSO}      & \textbf{RPSO}     & \textbf{PSO}      & \textbf{RPSO}     & \textbf{PSO}       & \textbf{RPSO}     \\
\hline
\textbf{Ackley1}                                       & 0.0248            & \textbf{0.0145}   & 0.0425            & \textbf{0.0102}   & 0.0152             & \textbf{0.0025}   \\
\textbf{Alpine1}                                       & 0.0013            & \textbf{0.0001}   & 0.0001            & \textbf{0.0001}   & 0.0010             & \textbf{0.0001}   \\
\textbf{Brent}                                         & \textbf{0.0000}   & 0.0000            & \textbf{0.0000}   & 0.0000            & \textbf{0.0000}    & 0.0000     \\
\textbf{Brown}                                         & 0.0000            & \textbf{0.0000}   & 0.0000            & \textbf{0.0000}   & 0.0002             & \textbf{0.0000}   \\
\textbf{ChungReynolds}                                 & \textbf{0.0000}   & 0.0000            & \textbf{0.0000}   & 0.0000            & 0.0000             & \textbf{0.0000}   \\
\textbf{Csendes}                                       & \textbf{0.0000}   & \textbf{0.0000}   & \textbf{0.0000}   & \textbf{0.0000}   & \textbf{0.0000}    & \textbf{0.0000}   \\
\textbf{Deb1}                                          & 0.0004            & \textbf{0.0000}   & 0.0003            & \textbf{0.0000}   & 0.0001             & \textbf{0.0000}   \\
\textbf{Deb2}                                          & 0.0001            & \textbf{0.0000}   & 0.0003            & \textbf{0.0000}   & 0.0001             & \textbf{0.0000}   \\
\textbf{DixonPrice}                                    & 0.0004            & \textbf{0.0000}   & 0.0003            & \textbf{0.0000}   & 0.0001             & \textbf{0.0000}   \\
\textbf{DropWave}                                      & 0.0011            & \textbf{0.0000}   & 0.0046            & \textbf{0.0001}   & 0.0026             & \textbf{0.0000}   \\
\textbf{EggHolder}                                     & 0.8794            & \textbf{0.0000}   & 0.9294            & \textbf{0.0001}   & 0.8795             & \textbf{0.0000}   \\
\textbf{Exponential}                                   & \textbf{0.0000}   & \textbf{0.0000}   & \textbf{0.0000}   & \textbf{0.0000}   & \textbf{0.0000}    & \textbf{0.0000}   \\
\textbf{Giunta}                                        & 0.0000            & \textbf{0.0000}   & 0.0000            & \textbf{0.0000}   & \textbf{0.0000}    & 0.0000            \\
\textbf{Mishra1}                                       & \textbf{0.0000}   & 0.0000            & \textbf{0.0000}   & 0.0000            & \textbf{0.0000}    & 0.0000            \\
\textbf{Mishra2}                                       & \textbf{0.0000}   & 0.0011            & \textbf{0.0000}   & 0.0002            & \textbf{0.0000}    & 0.0001            \\
\textbf{Periodic}                                      & 0.0004            & \textbf{0.0000}   & \textbf{0.0000}   & 0.0000            & 0.0003             & \textbf{0.0000}   \\
\textbf{PowellSum}                                     & \textbf{0.0000}   & \textbf{0.0000}   & \textbf{0.0000}   & \textbf{0.0000}   & \textbf{0.0000}    & \textbf{0.0000}   \\
\textbf{Qing}                                          & 0.0004            & \textbf{0.0000}   & 0.0002            & \textbf{0.0000}   & 0.0001             & \textbf{0.0000}   \\
\textbf{Rastrigin}                                     & 0.2921            & \textbf{0.0004}   & 0.1303            & \textbf{0.0014}   & 0.4009             & \textbf{0.0020}   \\
\textbf{Rosenbrock}                                    & 0.0029            & \textbf{0.0001}   & 0.0005            & \textbf{0.0000}   & 0.0027             & \textbf{0.0001}   \\
\textbf{Salomon}                                       & 0.0147            & \textbf{0.0003}   & 0.0133            & \textbf{0.0007}   & 0.0107             & \textbf{0.0010}   \\
\textbf{SchumerSteiglitz}                              & 0.0000            & \textbf{0.0000}   & 0.0000            & \textbf{0.0000}   & 0.0000             & \textbf{0.0000}   \\
\textbf{Sphere}                                        & \textbf{0.0000}   & 0.0000            & 0.0002            & \textbf{0.0000}   & 0.0002             & \textbf{0.0000}   \\
\textbf{Step}                                          & \textbf{0.0000}   & \textbf{0.0000}   & \textbf{0.0000}   & \textbf{0.0000}   & \textbf{0.0000}    & \textbf{0.0000}   \\
\textbf{StepInt}                                       & \textbf{0.0000}   & \textbf{0.0000}   & \textbf{0.0000}   & \textbf{0.0000}   & \textbf{0.0000}    & \textbf{0.0000}   \\
\textbf{SumSquares}                                    & 0.0003            & \textbf{0.0000}   & 0.0002            & \textbf{0.0000}   & \textbf{0.0000}    & 0.0000            \\
\textbf{Trid}                                          & \textbf{0.0000}   & 0.0000            & 0.0000            & \textbf{0.0000}   & 0.0000             & \textbf{0.0000}   \\
\textbf{Vincent}                                       & 0.0001            & \textbf{0.0000}   & 0.0001            & \textbf{0.0000}   & 0.0001             & \textbf{0.0000}   \\
\textbf{WWavy}                                         & \textbf{0.0000}   & \textbf{0.0000}   & \textbf{0.0000}   & \textbf{0.0000}   & \textbf{0.0000}    & \textbf{0.0000}   \\
\hline
\end{tabular}
}
\newline

\end{table}

Table \ref{tab:t3} presents the results of a comparative performance analysis of PSO and RPSO under different iteration limits together with other stopping criteria. The findings reveal that as the number of iterations increases, the performance advantage of RPSO over PSO becomes more pronounced. This trend suggests that the integration of the radial epiderivative mechanism allows RPSO to more effectively exploit additional iterations by guiding the particles in more informed and promising directions within the search space. As a result, RPSO is able to refine its solutions more efficiently over time, leading to enhanced overall performance, especially in scenarios that benefit from extended computational effort.

Table \ref{tab:t4} presents statistical data that compare the computational execution times of the PSO and RPSO algorithms. While RPSO generally requires slightly more computational time than PSO, this increase in runtime is modest and well justified by its improved performance in terms of solution quality. Specifically, RPSO consistently achieves smaller optimality gaps in various instances of problems, indicating more accurate and reliable solutions. This trade-off between runtime and solution accuracy highlights the practical advantage of RPSO, especially in applications where higher precision is prioritized over marginal differences in execution time.

\begin{table}
\caption{CPU time of PSO, RPSO, CC, and RCC in seconds}
\label{tab:t4}
\centering
\resizebox{\columnwidth}{!}{
\begin{tabular}{lcccc}
   \hline
\multicolumn{1}{c}{\multirow{2}{*}{\textbf{Problem}}} & \multicolumn{2}{c}{\textbf{Avg}} & \multicolumn{2}{c}{\textbf{Std Dev}} \\
\multicolumn{1}{c}{}                                   & \textbf{PSO}   & \textbf{RPSO}   & \textbf{PSO}     & \textbf{RPSO}     \\
\hline
\textbf{Ackley1}                                       & 0.008679       & 22.42503        & 0.005709         & 34.59871          \\
\textbf{Alpine1}                                       & 0.00639        & 16.22912        & 0.004006         & 22.18359          \\
\textbf{Brent}                                         & 0.013889       & 9.90954         & 0.008515         & 12.67603    \\
\textbf{Brown}                                         & 0.007881       & 5.784871        & 0.004768         & 8.943827          \\
\textbf{ChungReynolds}                                 & 0.006791       & 17.5923         & 0.004426         & 23.41651          \\
\textbf{Csendes}                                       & 0.00419        & 0.675348        & 0.002644         & 0.617709          \\
\textbf{Deb1}                                          & 0.006367       & 14.03301        & 0.00445          & 21.30635          \\
\textbf{Deb2}                                          & 0.006294       & 10.5939         & 0.00411          & 17.56088          \\
\textbf{DixonPrice}                                    & 0.008308       & 15.4246         & 0.005446         & 24.59655          \\
\textbf{DropWave}                                      & 0.009648       & 14.60314        & 0.006193         & 20.25077          \\
\textbf{EggHolder}                                     & 0.013651       & 28.00099        & 0.008268         & 51.63781          \\
\textbf{Exponential}                                   & 0.003543       & 0.31831         & 0.002309         & 0.303841          \\
\textbf{Giunta}                                        & 0.008163       & 9.170515        & 0.004716         & 13.68912          \\
\textbf{Mishra1}                                       & 0.005768       & 9.304859        & 0.003891         & 9.576343          \\
\textbf{Mishra2}                                       & 0.007377       & 22.3722         & 0.004186         & 29.95259          \\
\textbf{Periodic}                                      & 0.007218       & 14.84257        & 0.004937         & 25.20907          \\
\textbf{PowellSum}                                     & 0.004307       & 0.372636        & 0.003193         & 0.397355          \\
\textbf{Qing}                                          & 0.007696       & 68.496          & 0.004954         & 89.56419          \\
\textbf{Rastrigin}                                     & 0.007342       & 12.63646        & 0.005214         & 19.11386          \\
\textbf{Rosenbrock}                                    & 0.006766       & 13.27298        & 0.004747         & 21.56878          \\
\textbf{Salomon}                                       & 0.008379       & 32.16354        & 0.006535         & 36.30459          \\
\textbf{SchumerSteiglitz}                              & 0.006123       & 19.62367        & 0.003695         & 24.77052          \\
\textbf{Sphere}                                        & 0.00654        & 7.7528          & 0.004545         & 17.38545          \\
\textbf{Step}                                          & 0.003764       & 4.74353         & 0.002823         & 4.799111          \\
\textbf{StepInt}                                       & 0.007          & 1.264115        & 0.005155         & 1.673666          \\
\textbf{SumSquares}                                    & 0.00797        & 9.160162        & 0.005114         & 14.04453          \\
\textbf{Trid}                                          & 0.007483       & 11.13656        & 0.005224         & 20.01522          \\
\textbf{Vincent}                                       & 0.006169       & 9.305598        & 0.004115         & 12.79302          \\
\textbf{WWavy}                                         & 0.004097       & 0.834549        & 0.003363         & 1.282196          \\
   \hline
\end{tabular}
}

\end{table}

Table \ref{tab:t5} investigates the performance of the RPSO algorithm under different particle population sizes. The results demonstrate that RPSO maintains high solution quality even when the number of particles is relatively small, highlighting its efficiency in resource-constrained scenarios. This robustness is largely attributed to the integration of the radial epiderivative mechanism, which effectively guides the search process toward promising areas of the solution space. Additionally, in population-based metaheuristics, the exchange of information among particles plays a critical role in driving convergence. The radial epiderivative enhances this interaction by enabling more accurate and meaningful directional guidance, thereby improving the algorithm’s capacity to exploit useful information and converge more rapidly toward optimal or near-optimal solutions.

\begin{sidewaystable}
\caption{Computational results based on different particle sizes of PSO and RPSO}
\label{tab:t5}
\centering
\resizebox{\columnwidth}{!}{
\small
\begin{tabular}{lcccccccccccc}
\hline
\multicolumn{1}{c}{\multirow{2}{*}{\textbf{Problems}}} & \multicolumn{2}{c}{\textbf{\#ofParticles   = 5}} & \multicolumn{2}{c}{\textbf{\#ofParticles = 10}} & \multicolumn{2}{c}{\textbf{\#ofParticles = 20}} & \multicolumn{2}{c}{\textbf{\#ofParticles = 30}} & \multicolumn{2}{c}{\textbf{\#ofParticles = 50}} & \multicolumn{2}{c}{\textbf{\#ofParticles = 100}} \\
\multicolumn{1}{c}{}                                   & \textbf{PSO}           & \textbf{RPSO}          & \textbf{PSO}          & \textbf{RPSO}          & \textbf{PSO}          & \textbf{RPSO}          & \textbf{PSO}          & \textbf{RPSO}          & \textbf{PSO}          & \textbf{RPSO}          & \textbf{PSO}           & \textbf{RPSO}          \\
\hline
\textbf{Ackley1}                                       & 0.821815               & 0.038662               & 0.717003              & 0.071462               & 0.042503              & 0.030335               & 0.064616              & 0.026659               & 0.015244              & 0.014452               & 0.024832               & 0.002486               \\
\textbf{Alpine1}                                       & 0.007553               & 8.58E-05               & 0.002606              & 0.00025                & 0.000355              & 0.000127               & 0.001044              & 0.000322               & 8.41E-05              & 7.59E-05               & 0.001312               & 6.99E-05               \\
\textbf{Brent}                                         & 1.38E-87               & 0.001415               & 1.38E-87              & 9.15E-05               & 1.38E-87              & 2.82E-05               & 1.38E-87              & 1.75E-06               & 1.38E-87              & 2.72E-05               & 1.38E-87               & 3.16E-07              \\
\textbf{Brown}                                         & 0.002464               & 0.000614               & 0.002441              & 4.61E-05               & 0.00068               & 1.51E-06               & 0.000955              & 1.08E-05               & 3.79E-05              & 5.45E-05               & 1.09E-05               & 0.000204               \\
\textbf{ChungReynolds}                                 & 4.71E-05               & 4.8E-05                & 2.1E-07               & 3.54E-06               & 1.98E-12              & 1.29E-10               & 1.31E-10              & 4.12E-10               & 2.62E-08              & 7.74E-08               & 5.61E-12               & 1.48E-11               \\
\textbf{Csendes}                                       & 0                      & 0                      & 0                     & 0                      & 0                     & 0                      & 0                     & 0                      & 0                     & 0                      & 0                      & 0                      \\
\textbf{Deb1}                                          & -0.99913               & -0.99991               & -0.99865              & -1                     & -0.99883              & -1                     & -0.99947              & -0.99999               & -0.9993               & -1                     & -0.9997                & -1                     \\
\textbf{Deb2}                                          & -0.97428               & -0.99993               & -0.99586              & -0.99999               & -0.99974              & -1                     & -0.99916              & -0.99999               & -0.99962              & -1                     & -0.99979               & -0.99999               \\
\textbf{DixonPrice}                                    & 0.007304               & 0.002034               & 0.002071              & 0.000263               & 0.001729              & 0.000164               & 0.000286              & 7.66E-06               & 5.68E-05              & 1.57E-05               & 0.00042                & 5.85E-06               \\
\textbf{DropWave}                                      & -0.93624               & -0.98156               & -0.93623              & -0.99557               & -0.98064              & -0.99943               & -0.98824              & -0.99994               & -0.9978               & -1                     & -0.99618               & -1                     \\
\textbf{EggHolder}                                     & -114.891               & -951.149               & -66.8409              & -932.735               & -66.8425              & -959.631               & -66.8435              & -956.895               & -66.8416              & -959.641               & -66.8434               & -959.632               \\
\textbf{Exponential}                                   & -1                     & -1                     & -1                    & -1                     & -1                    & -1                     & -1                    & -1                     & -1                    & -1                     & -1                     & -1                     \\
\textbf{Giunta}                                        & 0.064734               & 0.064741               & 0.064606              & 0.064521               & 0.064672              & 0.064514               & 0.064472              & 0.064472               & 0.064483              & 0.064471               & 0.064474               & 0.064472               \\
\textbf{Mishra1}                                       & 2                      & 2.000001               & 2                     & 2                      & 2                     & 2                      & 2                     & 2                      & 2                     & 2                      & 2                      & 2                      \\
\textbf{Mishra2}                                       & 2                      & 2.014545               & 2                     & 2.006724               & 2                     & 2.001976               & 2                     & 2.000442               & 2                     & 2.002531               & 2                      & 2.000619               \\
\textbf{Periodic}                                      & 0.908089               & 0.901359               & 0.903848              & 0.900165               & 0.900717              & 0.900006               & 0.900739              & 0.900133               & 0.900543              & 0.9                    & 0.9                    & 0.900004               \\
\textbf{PowellSum}                                     & 0                      & 0                      & 0                     & 0                      & 0                     & 0                      & 0                     & 0                      & 0                     & 0                      & 0                      & 0                      \\
\textbf{Qing}                                          & 0.075714               & 0.003295               & 0.007726              & 0.000709               & 0.00354               & 1.4E-06                & 0.000914              & 8.61E-05               & 0.003138              & 1.43E-05               & 0.000122               & 2.18E-06               \\
\textbf{Rastrigin}                                     & 1.127427               & 0.01777                & 1.84131               & 0.033036               & 0.602677              & 0.001768               & 0.149944              & 0.002201               & 0.292112              & 0.002036               & 0.130333               & 0.000364               \\
\textbf{Rosenbrock}                                    & 0.002924               & 0.001776               & 0.002908              & 5.74E-05               & 0.002966              & 4.76E-05               & 0.001481              & 0.000273               & 0.000501              & 0.000363               & 0.002694               & 4.49E-05               \\
\textbf{Salomon}                                       & 0.100527               & 0.004902               & 0.014678              & 0.016959               & 0.013301              & 0.000659               & 0.014035              & 0.001062               & 0.021703              & 0.000343               & 0.010694               & 0.000355               \\
\textbf{SchumerSteiglitz}                              & 7.45E-06               & 1.37E-06               & 3.01E-08              & 2.87E-07               & 8.57E-07              & 1.68E-09               & 4.51E-07              & 1.74E-09               & 5.18E-10              & 6.27E-10               & 1.73E-09               & 3.18E-11               \\
\textbf{Sphere}                                        & 0.005374               & 0.001326               & 0.002034              & 0.001204               & 0.000141              & 8.35E-05               & 0.000139              & 2.07E-05               & 3.24E-06              & 1.32E-05               & 1.51E-05               & 5.49E-05               \\
\textbf{Step}                                          & 0                      & 0                      & 0                     & 0                      & 0                     & 0                      & 0                     & 0                      & 0                     & 0                      & 0                      & 0                      \\
\textbf{StepInt}                                       & 13                     & 13                     & 13                    & 13                     & 13                    & 13                     & 13                    & 13                     & 13                    & 13                     & 13                     & 13                     \\
\textbf{SumSquares}                                    & 0.002428               & 0.006272               & 0.002364              & 0.013582               & 2.25E-06              & 0.000384               & 0.000845              & 4.04E-05               & 0.000104              & 5.5E-05                & 0.000171               & 7.26E-06               \\
\textbf{Trid}                                          & -1.99769               & -1.99994               & -1.99901              & -1.99981               & -1.99984              & -1.99993               & -1.9999               & -1.99998               & -1.99996              & -1.99998               & -1.99997               & -1.99999               \\
\textbf{Vincent}                                       & -1.99745               & -1.99984               & -1.99954              & -2                     & -1.99985              & -2                     & -1.99844              & -2                     & -1.9997               & -2                     & -1.99982               & -2                     \\
\textbf{WWavy}                                         & 0                      & 0                      & 0                     & 0                      & 0                     & 0                      & 0                     & 0                      & 0                     & 0                      & 0                      & 0                      \\

\hline
\end{tabular}
}
\newline
\end{sidewaystable}
Table \ref{tab:t6} provides a detailed comparative analysis of the performance of PSO and RPSO across 216 experiments conducted under various parameter configurations. To ensure a more meaningful comparison, the table excludes the experiments in which both algorithms produced identical best results—either in terms of absolute objective values or percentage differences. The remaining data offer strong evidence that RPSO consistently outperforms PSO in terms of convergence quality. For example, for Alpine 1, RPSO prevailed in 194 out of 216 experiments with different parameters, while PSO prevailed in only 22 of these experiments. In this table, in a row that does not total 216, it is understood that no one method outperformed the other in the relevant experiments. The percentages in this table show the proportion of experiments in which one method prevailed. Specifically, RPSO not only finds better solutions more frequently but also exhibits greater robustness across diverse parameter settings, highlighting its enhanced capability in navigating complex search spaces.

The concave functions included in the problem set were solved using the SciPy optimization module \citep{2020SciPy-NMeth}, which served as a benchmark for evaluating solution quality. The results obtained from SciPy were then compared with those generated by the PSO and RPSO algorithms. As illustrated in Table \ref{tab:txy}, RPSO consistently delivers solutions that are significantly closer to those produced by SciPy, indicating a strong alignment with optimal or near-optimal values. In contrast, standard PSO exhibits larger deviations from the benchmark results. These findings highlight RPSO's superior ability to handle concave optimization problems, demonstrating its robustness and improved convergence behavior over traditional PSO.

\begin{table}[h!]
\centering
\caption{Comparisons on PSO-RPSO experiments with different parameters}
\label{tab:t6}
\renewcommand{\arraystretch}{0.9} 
\setlength{\tabcolsep}{6pt} 
\small 
\begin{tabular}{lcccc}
\toprule
\textbf{Problem} & \textbf{RPSO} & \textbf{PSO} & \textbf{RPSO\%} & \textbf{PSO\%} \\
\midrule
\textbf{Ackley1}          & 125  & 91   & 57.87 & 42.13 \\
\textbf{Alpine1}          & 194  & 22   & 89.81 & 10.19 \\
\textbf{Brent}            & 28   & 188  & 12.96 & 87.04 \\
\textbf{Brown}            & 105  & 111  & 48.61 & 51.39 \\
\textbf{ChungReynolds}    & 119  & 97   & 55.09 & 44.91 \\
\textbf{Csendes}          & 0    & 0    & 0.00  & 0.00  \\
\textbf{Deb1}             & 206  & 10   & 95.37 & 4.63  \\
\textbf{Deb2}             & 189  & 27   & 87.50 & 12.50 \\
\textbf{DixonPrice}       & 131  & 85   & 60.65 & 39.35 \\
\textbf{DropWave}         & 181  & 35   & 83.80 & 16.20 \\
\textbf{EggHolder}        & 216  & 0    & 100.00 & 0.00 \\
\textbf{Exponential}      & 0    & 1    & 0.00  & 100.00 \\
\textbf{Giunta}           & 82   & 134  & 37.96 & 62.04 \\
\textbf{Mishra1}          & 4    & 212  & 1.85  & 98.15 \\
\textbf{Mishra2}          & 6    & 210  & 2.78  & 97.22 \\
\textbf{Periodic}         & 197  & 19   & 91.20 & 8.80  \\
\textbf{PowellSum}        & 0    & 0    & 0.00  & 0.00  \\
\textbf{Qing}             & 145  & 71   & 67.13 & 32.87 \\
\textbf{Rastrigin}        & 191  & 25   & 88.43 & 11.57 \\
\textbf{Rosenbrock}       & 137  & 79   & 63.43 & 36.57 \\
\textbf{Salomon}          & 197  & 19   & 91.20 & 8.80  \\
\textbf{SchumerSteiglitz} & 105  & 111  & 48.61 & 51.39 \\
\textbf{Sphere}           & 101  & 115  & 46.76 & 53.24 \\
\textbf{Step}             & 56   & 3    & 94.92 & 5.08  \\
\textbf{StepInt}          & 108  & 58   & 65.06 & 34.94 \\
\textbf{SumSquares}       & 101  & 115  & 46.76 & 53.24 \\
\textbf{Trid}             & 77   & 139  & 35.65 & 64.35 \\
\textbf{Vincent}          & 194  & 22   & 89.81 & 10.19 \\

\textbf{WWavy}            & 25   & 23   & 52.08 & 47.92 \\
\bottomrule
\end{tabular}
\end{table}

\begin{sidewaystable}
\caption{Concave Results Comparison}
\label{tab:txy}

\centering
\resizebox{\columnwidth}{!}{
\small
\begin{tabular}{lrrrrr}
\toprule
\textbf{Problem} & \textbf{PSO} & \textbf{RPSO} & \textbf{Scipy} & \textbf{Gap-PSO} & \textbf{Gap-RPSO} \\
\midrule
\textbf{Brown}            & $-5.90295 \times 10^{20}$ & $-5.90295 \times 10^{20}$ & $-5.90295 \times 10^{20}$ & \textbf{0}   & \textbf{0}   \\
\textbf{Brent}            & -430.9689                & -799.9913                & -800                      & 0.4607                         & \textbf{$1.091 \times 10^{-5}$} \\
\textbf{ChungReynolds}    & -2289.031                & $-3.99999 \times 10^{8}$ & $-4.00000 \times 10^{8}$  & 0.99999                       & \textbf{$2.199 \times 10^{-6}$} \\
\textbf{DixonPrice}       & -5687.8225               & -44220.2004              & -44221                    & 0.8714                        & \textbf{$1.808 \times 10^{-5}$} \\
\textbf{Exponential}      & 0.5668                   & 0.3684                   & 0.3679                    & 0.1454                        & \textbf{0.0004} \\
\textbf{PowellSum}        & -1.3796                  & -1.9981                  & -2                        & 0.2068                        & \textbf{0.0006} \\
\textbf{SchumerSteiglitz} & -952.7147                & $-1.99999 \times 10^{8}$ & $-2.00000 \times 10^{8}$  & 0.99999                       & \textbf{$3.054 \times 10^{-7}$} \\
\textbf{Sphere}           & -52.4288                 & -52.4054                 & -52.4288                  & \textbf{0}                    & 0.0004 \\
\textbf{SumSquares}       & -61.3821                 & -299.9721                & -300                      & 0.7928                        & \textbf{$9.281 \times 10^{-5}$} \\
\bottomrule
\end{tabular}
}
\end{sidewaystable}

\section{Conclusion}
This paper proposes a generic solution method for box constrained nonconvex optimization problems which uses the descent direction and the global optimality criteria provided by the radial epiderivatives of objective functions. Probably these features of the radial epiderivative are used in the presented method for the first time. In this paper, we combined the above-mentioned abilities of the radial epideraivative with two search methods known as CC and PSO. The obtained results demonstrate that the hybrid methods (RCC and RPSO) derived by this way, outperform the sole search methods. This encourages us to test in feature new combinations for the radial epiderivative based method with different search algorithms and also to apply the radial epiderivative based search method to new classes of nonsmooth and nonconvex problems with objective functions whose radial epiderivatives can be calculated more accurately. We also believe that the methods developed in this paper, can be used to generate new hybrid exact solution methods in constrained optimization by combining them with augmented Lagrangian based primal-dual methods (see e.g. \citep{AlpaslanK2020,BagirovOK2019,BulbulK2025,BulbulK2021,Gasimov2002,GasimovU2007,KasimbeyliUR2009}).

\newpage


\section{Acknowledgements}
This study was supported by Eskisehir Technical University Scientific Research Projects Commission under grants no: 24ADP030 and 24ADP031.
\newpage
\appendix
\section*{Test Problems}

\textbf{Ackley1}:
$$f(X)=-20 \ e^{\displaystyle\left(-0.2\times\sqrt{\frac{1}{n}\sum_{i=1}^{n}x^2_i}\right)}-e^{\displaystyle\left[\frac{1}{n}\sum_{i=1}^{n}\cos(2\pi x_i)\right]}+20+e^{(1)}$$ 
\textbf{Alpine1}:
$$f(X)=\sum_{i=1}^{n}\left|x_i\sin\left(x_i\right)+0.1x_i\right|$$
\textbf{Brown}:
$$f(X)=\sum^{n-1}_{i=1}\Big(x^2_i\Big)^{\displaystyle \left(x^2_{i+1}+1\right)}+\Big(x^2_{i+1}\Big)^{\displaystyle\left(x^2_{i}+1\right)}$$
\textbf{ChungReynolds}:
$$f(X)=\left(\sum_{i=1}^{n}x^2_i\right)^2$$
\textbf{Csendes}:
$$f(X)=\begin{cases}  \sum_{i=1}^{n}x^6_i\left[\sin\left(\frac{1}{x_i}\right)+2\right] & \text{ if } \prod^{n}_{i=1} x_i \neq 0 \\  0 & \text{ otherwise }  \end{cases}$$
\textbf{Deb1}:
$$f(X)=-\frac{1}{n}\sum^n_{i=1}\sin^6\left(5 \pi x_i\right)$$
\textbf{Deb2}:
$$f(X)=-\frac{1}{n}\sum^n_{i=1}\sin^6\left[5 \pi \left(x^{3/4}_i-0.05\right)\right]$$
\textbf{DixonPrice}:
$$f(X)=\left(x_1-1\right)^2+\sum^n_{i=2}i\left(2x^2_i-x_{i-1}\right)^2$$
\textbf{DropWave}:
$$f(X)=-\frac{1+\cos\displaystyle\left(12\sqrt{x^2_1+x^2_2}\right)}{\displaystyle\frac{1}{2}\left(x^2_1+x^2_2\right)+2}$$
\textbf{EggHolder}:
$$f(X)=\sum_{i=1}^{n-1}\left[-x_i\sin\left(\sqrt{\left|x_i-x_{i+1}-47\right|}\right)-\left(x_{i+1}+47\right)\sin\left(\sqrt{\left|0.5x_i+x_{i+1}+47\right|}\right)\right]$$
\textbf{Exponential}:
$$f(X)=\text{exp}\left(-0.5\sum^{n}_{i=1} x^2_i\right)$$
\textbf{Giunta}:
$$f(X)=0.6+\sum_{i=1}^{n}\left[\sin^2\left(1-\frac{16}{15}x_i\right)-\frac{1}{50}\sin\left(4-\frac{64}{15}x_i\right)-\sin\left(1-\frac{16}{15}x_i\right)\right]$$
\textbf{Mishra1}:
$$f(X)=\left(1+g_n\right)^{g_n}$$
where: $$g_n=n-\sum_{i=1}^{n-1}x_i$$
\textbf{Mishra2}:
$$f(X)=\left(1+g_n\right)^{g_n}$$
where: $$g_n=n-\sum_{i=1}^{n-1}\frac{\displaystyle \left(x_i+x_{i+1}\right)}{\displaystyle 2}$$
\textbf{Periodic}:
$$f(X)=1+\sin^2\left(x_1\right)+\sin^2\left(x_2\right)-0.1e^{\left(-x_1^2-x_2^2\right)}$$
\textbf{PowellSum}:
$$f(\mathbf{x}) = \sum_{i=1}^n |\mathbf{x}_i|^{i+1}$$
\textbf{Qing}:
$$f(X)=\sum_{i=1}^{n}\left(x^2_i-i\right)^2$$
\textbf{Rastrigin}:
$$f(X)=\sum_{i=1}^{n}\left[x^2_i-10\cos\left(2\pi x_i\right)+10\right]$$
\textbf{Rosenbrock}:
$$f(X)=\sum^{n-1}_{i=1}\left[100 \left|x_{i+1} - x_i^2\right| + \left(1-x_i\right)^2\right]$$
\textbf{Salomon}:
$$f(X)=1-\cos\left(2\pi \left \| x \right \|\right)+0.1\left \| x \right \|$$
\textbf{SchumerSteiglitz}:
$$f(X)=\sum^n_{i=1} x^4_i$$
\textbf{Sphere}:
$$f\left( {X} \right) = \sum\limits_{i = 1}^n {{x_i}^2}$$
\textbf{Step}:
$$f(X)=\sum^n_{i=1}\left \lfloor \left|x_i\right| \right \rfloor$$
\textbf{StepInt}:
$$f(X)=25+\sum^n_{i=1}\left \lfloor x_i \right \rfloor$$
\textbf{SumSquares}:
$$f(X)=\sum_{i=1}^{n}ix^2_i$$
\textbf{Trid}:
$$f(X)=\sum_{i=1}^{n}\left(x_i-1\right)^2-\sum_{i=2}^{n}x_i \cdot x_{i-1}$$
\textbf{Vincent}:
$$f(X)=-\frac{\displaystyle 1}{\displaystyle n}\sum^n_{i=1} \sin\left[10 \log\left(x_i\right)\right]$$
\textbf{WWavy}:
$$f(X)=\frac{\displaystyle 1}{\displaystyle n}\sum^n_{i=1}1-\cos\left(kx_i\right)e^{-\frac{1}{2}x_i^2}$$
where: $$k=10$$
\textbf{Brent}:
$$f(X)=\left(x_1+10\right)^2+\left(x_2+10\right)^2+e^{-x^2_1-x^2_2}$$

\bibliographystyle{elsarticle-harv} 
\bibliography{cas-refs}






\end{document}